\documentclass[11pt]{article}
\usepackage{url} %this package should fix any errors with URLs in refs.
\usepackage{lineno}
\usepackage[margin=2cm]{geometry}
\usepackage{authblk}
\usepackage{soul}
\usepackage{amsmath}
\usepackage{algorithm}
\usepackage{algpseudocode}
\usepackage{caption}
\usepackage{subcaption}
\usepackage{amssymb}
\usepackage{graphicx}
\usepackage{biblatex}
\bibliography{references.bib}

\begin{document}
\title{Effective Preconditioners for Mixed-Dimensional Scalar Elliptic Problems}
\author[1,2]{Xiaozhe Hu}
\author[2]{Eirik Keilegavlen\thanks{e-mail: Eirik.Keilegavlen@uib.no},}
\author[2]{Jan M.~Nordbotten}
\affil[1]{Department of Mathematics, Tufts University, Medford, MA 02155, USA}
\affil[2]{Center for Modeling of Coupled Subsurface Dynamics, Department of Mathematics, University of Bergen, 5020, Bergen, Norway}
\date{June 2022}
\maketitle

\begin{abstract}
Discretization of flow in fractured porous media commonly lead to large systems of linear  equations that require dedicated solvers.
In this work, we develop an efficient linear solver and its practical implementation for mixed-dimensional scalar elliptic problems. We design an effective preconditioner based on approximate block factorization and algebraic multigrid techniques. Numerical results on benchmarks with complex fracture structures demonstrate the effectiveness of the proposed linear solver and its robustness with respect to different physical and discretization parameters.
\end{abstract}

%% ------------------------------------------------------------------------ %%
%
%  TEXT
%
%% ------------------------------------------------------------------------ %%

\section{Introduction}
Mixed-dimensional scalar elliptic equations are the backbone of models across many applications wherein potentially intersecting thin layers are embedded into a material. This application is of particular relevance in subsurface flows, where the thin layers correspond to fracture networks embedded in a porous material. 

Flow in fractured porous media can be modeled in several ways. However, three of the main modeling approaches result in essentially mixed-dimensional elliptic equations. In sequence, these are 1) Direct mixed-dimensional modeling \cite{frih2012modeling,boon2018robust,boon2020functional}, 2) Discrete fracture networks \cite{erhel2009,hyman2015dfnworks}, and 3) Low-order numerical methods for equidimensional fracture models \cite{karimi2004efficient,sandven2012efficient}. These modeling approaches are presented in detail in the recent review paper \cite{berre2019flow}, and two community benchmark studies have recently been conducted \cite{Flemisch2016a,berre2021verification}. The relative merits of the three modeling approaches have been extensively analyzed in the literature. However, for our purpose, it suffices to point out that after discretization by standard numerical methods, the resulting algebraic linear systems have a similar structure. 

Despite the importance of modeling flow in fractured porous media, the computational cost of linear solvers has not received much attention, although linear solvers for discrete fracture networks recently were studied in \cite{greer2022comparison}.
Indeed, in recent benchmark studies, the computational cost is simply represented by proxy in terms of the condition number of the system matrix. The construction, or even existence, of efficient iterative solvers can not be taken for granted. The material contrasts between porous rock and fracture may be significant, which is known to cause a challenge for many iterative algorithms, see~\cite{chan2000,mandel1996,nepomnyaschikh1999,oswald1999,graham1999,trottenberg2000multigrid} and references therein.  

This technical note presents an efficient iterative linear solver for the prototypical mixed-dimensional scalar elliptic problem, which serves as a building block for solving other mixed-dimensional problems, e.g., \cite{budivsa2021block} and~\cite{budisa2020mixed}. We explore the block structure of the results linear systems and design effective block preconditioners based on approximate block factorization and algebraic multigrid (AMG) methods. The successful application requires several essential steps. Firstly, because the blocks corresponding to fractures are usually smaller in size and lower in dimension, we build an approximation Schur complement on the domain via a simple diagonal approximation to the fracture block. Secondly, for the diagonal blocks in the block preconditioners, we employ existing AMG methods to approximately invert the diagonal blocks via several steps of V-cycles. Finally, the block preconditioners are used to accelerate the Krylov iterative methods, which results in an effective and robust linear solver. %We present numerical results on several benchmark problems and report the performance of the proposed linear solver to demonstrate its effectiveness numerically. 

This combination of stable discretization and efficient linear solvers has not been reported in the literature to the best of our knowledge. The novelty of this technical note is identifying an efficient and robust iterative solver for mixed-dimensional scalar elliptic equations. We support this claim by applying the solver to one 2D benchmark problem with complex fracture networks~\cite{Flemisch2016a}, and to two 3D benchmark problems~\cite{berre2021verification} with structured and complex fracture networks, respectively, and show robustness (in terms of iteration count until convergence) across a range of material parameters and grid resolutions.

\section{Preliminaries}

\subsection{Mixed-dimensional Scalar Elliptic Problems}
\begin{figure}
    \centering
    \includegraphics[width=\textwidth]{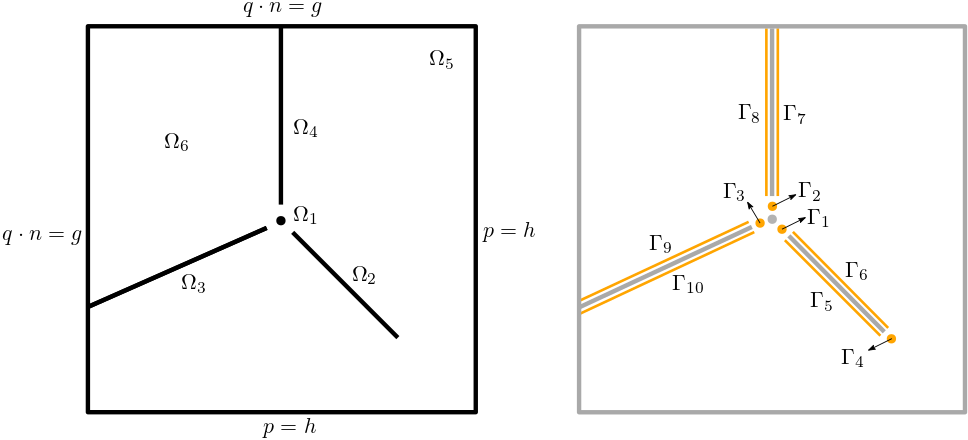}
    \caption{Left: A mixed-dimensional geometry consisting of three fractures, one intersection point, and with the top-dimensional domain split into two. Possible boundary conditions are indicated. Right: Geometry of the interfaces; $\Gamma_4$ is included as a placeholder for a no-flux condition on the tip of $\Omega_2$. }
    \label{fig:GeometryNotation}
\end{figure}

Our model for flow in fractured porous media is the same as used in the benchmark for 3d flow \cite{berre2021verification} and several other works, e.g., \cite{boon2018robust,glaser2022comparison}. 
We model the fractures as $(N-1)$-dimensional objects embedded in the N-dimensional host medium $\Omega\subset \mathbb{R}^N$, with fracture intersections forming lines and points of dimensions $(N-2)$ and $(N-3)$.
We refer to these geometric objects as subdomains and let $\Omega_i$ denote an arbitrary subdomain, so that the computational domain $\Omega=\bigcup_{i=1}^m \Omega_i$.
The subdomains are connected via interfaces, denoted $\Gamma_j, j=\{1,...,M\}$. If a $D$-dimensional interface $\Gamma_j$ connects subdomains $\Omega_{\hat{j}}$ and $\Omega_{\check{j}}$, of dimension $D+1$ and $D$, respectively, we have that geometrically $\Gamma_j=\Omega_{\check{j}}$, and furthermore denote the part of the boundary of $\Omega_{\hat{j}}$ coinciding with  $\Gamma_j$ as $\Gamma_j=\partial_j\Omega_{\hat{j}}$.
Finally, for $\Omega_i$, let $\hat{S}_i$ 
be the set 
of interfaces towards neighboring subdomains of higher 
dimension, such that $j\in \hat{S}_i$ if and only if $\check{j}=i$.
The geometry is illustrated in Figure \ref{fig:GeometryNotation}, where we note that a combination of fracture geometry and boundary condition can lead to subdomains (of dimension $N$ or lower) being disconnected from Dirichlet boundary conditions.

With this notion of geometry, the model for single-phase fluid flow in $\Omega_i$ is given in strong form as,
\begin{equation} \label{eqn:scalar-elliptic-domain}
    \nabla_{\parallel}\cdot(-K_{i,\parallel}\nabla_{\parallel}p_i) - \sum\limits_{j\in \hat{S}_i} \lambda_j=f_i .
\end{equation}
Here, $K_i$, $p_i$, and $f_i$ represent permeability, pressure, and source terms in subdomain $\Omega_i$, and subscript $\parallel$ indicates that operators are tangential to $\Omega_i$. 
The flux over $\Gamma_j$, denoted $\lambda_j$, is governed by a Darcy-type law on the form
\begin{equation} \label{eqn:scalar-elliptic-interface}
    \lambda_j = -\kappa_j (p_i - p_{\hat{j}}),
\end{equation}
The constant $\kappa_j$ can be considered related to the normal transmissivity across $\Gamma_j$.
Thus the terms $\lambda_j$ represent flow through $\Gamma_j$ from $\Omega_{\hat{j}}$ to $\Omega_i$.
This same flow through $\Gamma_j$ then appears as a Neumann boundary term on $\partial_j\Omega_{\hat{j}}$, 
\begin{equation}
    \label{eqn:jan-added-this}
    (-K_{\hat{j},\parallel}\nabla_{\parallel}p_ {\hat{j}})\cdot\mathbf{n}_{\hat{j}}= \lambda_j
\end{equation}
where $\mathbf{n}_{\hat{j}}$ is the outer normal vector on the boundary of $\Omega_{\hat{j}}$. 

Equations \eqref{eqn:scalar-elliptic-domain} remain valid also for the highest-dimensional domains (when $\hat{S}_i=\emptyset$, and as such the summation is void), and also for the lowest-dimensional point intersections (when there is no parallel directions, and as such the differential operators are void). 

\subsection{Discretization} 
We discretize the flow model following the unified framework developed in \cite{nordbotten2019unified}, wherein the subdomains are discretized independently with the interface fluxes treated as interior boundary terms. The framework is compatible with various discretization methods for elliptic equations; herein, we apply the finite volume multipoint flux approximation method \cite{aavatsmark2002introduction} to discretize subdomain problems, but mixed and virtual finite elements have also been used \cite{nordbotten2019unified}.

\section{Preconditioners for Mixed-dimensional Problems} 
In this section, we introduce the multigrid-based solver for solving the discretized mixed-dimensional problem~\eqref{eqn:scalar-elliptic-domain} and~\eqref{eqn:scalar-elliptic-interface}. Our solver takes advantage of the block structure of the linear systems after discretization.  More precisely, we use the block structure to perform a block factorization and develop preconditioners based on the block factorization via appropriate approximated Schur complement.  Finally, the multigrid method is used to obtain an efficient and practical preconditioner. 

\subsection{Block Structure of the Linear Systems of Equations}
Since \eqref{eqn:scalar-elliptic-domain} is defined in the domain and~\eqref{eqn:scalar-elliptic-interface} is defined on the interface, our the discretized linear systems naturally inherits the following two-by-two block structure,
\begin{equation} \label{eqn:2by2-system}
\mathcal{A} p = f  \quad \Leftrightarrow 
\begin{pmatrix}
A_{\Omega \Omega} & A_{\Omega\Gamma} \\
A_{\Gamma\Omega} & A_{\Gamma\Gamma}\\
\end{pmatrix}	
\begin{pmatrix}
p_{\Omega} \\
p_{\Gamma}
\end{pmatrix}
= 
\begin{pmatrix}
f_{\Omega} \\
f_{\Gamma}
\end{pmatrix},
\end{equation}
where the subscripts $\Omega$ and $\Gamma$ denote the blocks related to the domain and the interface, respectively. Correspondingly
$p_{\Omega}$ and $p_{\Gamma}$ denote the unknowns in all the subdomain $\Omega_i$ and all the interfaces $\Gamma_i$, respectively. The off-diagonal blocks $A_{\Omega\Gamma}$ and $A_{\Gamma\Omega}$ denote the interaction and coupling between the domain and the interfaces. As our implementation follows the unified discretization approach for fracture \cite{nordbotten2019unified}, we are guaranteed that $A_{\Omega\Gamma} = A_{\Gamma\Omega}^T$.  However, for the sake of generality, we will not exploit this property in the following discussion.  

A block factorization-based approach is a natural choice for solving the linear systems with block structure, such as~(\ref{eqn:2by2-system}).  There are two types of block factorization for a two-by-two block system which give different Schur complements.  Here, since $A_{\Gamma\Gamma}$ is defined on lower dimension interfaces, it is relatively smaller and easier to invert compared with $A_{\Omega\Omega}$.  Hence, we use the following form of block factorization.
\begin{equation}\label{eqn:block-fac-form}
\mathcal{A} = 
\begin{pmatrix}
	A_{\Omega \Omega} & A_{\Omega\Gamma} \\
	A_{\Gamma\Omega} & A_{\Gamma\Gamma}\\
\end{pmatrix}
= 
\begin{pmatrix}
I & A_{\Omega\Gamma}A_{\Gamma\Gamma}^{-1}  \\
0  & I	
\end{pmatrix}
\begin{pmatrix}
S_{\Omega\Omega} & 0 \\
0  & A_{\Gamma\Gamma}
\end{pmatrix}
\begin{pmatrix}
	I &  0  \\
	A_{\Gamma\Gamma}^{-1}A_{\Gamma\Omega}  & I	
\end{pmatrix} 
=: \mathcal{U} \mathcal{D} \mathcal{L},
\end{equation}
where $S_{\Omega\Omega}$ is the well-known Schur complement defined as
\begin{equation}\label{eqn:Schur-comp}
S_{\Omega\Omega} := A_{\Omega \Omega} - A_{\Omega \Gamma} A_{\Gamma \Gamma}^{-1} A_{\Gamma \Omega}.
\end{equation} 
Such a block factorization~(\ref{eqn:block-fac-form}) serves the building block of our proposed block preconditioner. 

\subsection{Factorization-based Block Preconditioner}
Based on the block factorization~\eqref{eqn:block-fac-form}, we can immediately propose several block preconditioners, e.g., $\mathcal{B}_D = \mathcal{D}^{-1}$ can be used as a block diagonal preconditioner, $\mathcal{B}_U := (\mathcal{UD})^{-1}$ and $\mathcal{B}_L := (\mathcal{DL})^{-1}$ can be used as block upper and lower triangular preconditioners, respectively.  For the sake of simplicity, in this subsection, we focus our discussion on $\mathcal{B}_L$ (and its variants).
% and would like to point out that similar discussion can be done for both $\mathcal{B}_D$ and $\mathcal{B}_U$ as well.

If we use $\mathcal{B}_L$ as a right preconditioner for Krylov iterative methods, such as general minimal residual (GMRes) method, we need to look at the spectrum of the preconditioned linear system $\mathcal{A}\mathcal{B}_L$.  From~\eqref{eqn:block-fac-form}, we have
%\begin{equation*}
$\mathcal{A}\mathcal{B}_L = \mathcal{UDL} (\mathcal{DL})^{-1} = \mathcal{U}$.
%\end{equation*}
This immediately implies that the eigenvalues of $\mathcal{A}\mathcal{B}_L$ are
%\begin{equation*}
$	\lambda(\mathcal{A}\mathcal{B}_L) = \lambda(\mathcal{U}) = 1$.
%\end{equation*}
Therefore, $\mathcal{B}_L$ is an efficient and robust preconditioner for solving~\eqref{eqn:2by2-system}.  
 
Note that the action of
\begin{equation*}
\mathcal{B}_L = 
\begin{pmatrix}
	S_{\Omega\Omega} &  0  \\
	A_{\Gamma\Omega}  & A_{\Gamma\Gamma}	
\end{pmatrix}^{-1}
\end{equation*}
requires computing the Schur complement $S_{\Omega\Omega}$, the inverse of the Schur complement $S_{\Omega\Omega}^{-1}$, and $A_{\Gamma\Gamma}^{-1}$.  All these steps could be quite expensive and we instead approximate them one by one.

For approximating the Schur complement, we use the fact that $A_{\Gamma\Gamma}$ in general is diagonally dominated and replace $A_{\Gamma\Gamma}$ in~\eqref{eqn:Schur-comp} by its diagonal $\operatorname{diag}(A_{\Gamma\Gamma})$.  Therefore, the approximate Schur complement is defined as
\begin{equation}\label{eqn:approx-Schur-comp}
	\widetilde{S}_{\Omega\Omega} := A_{\Omega \Omega} - A_{\Omega \Gamma} \left(\operatorname{diag}(A_{\Gamma \Gamma})\right)^{-1} A_{\Gamma \Omega}.
\end{equation} 
This approximation  leads to the following approximate block lower triangular preconditioner
\begin{equation*}
	\widetilde{\mathcal{B}}_L = 
	\begin{pmatrix}
		\widetilde{S}_{\Omega\Omega} &  0  \\
		A_{\Gamma\Omega}  & A_{\Gamma\Gamma}	
	\end{pmatrix}^{-1}.
\end{equation*}
Unfortunately, exactly inverting $\widetilde{\mathcal{B}}_L$ is still expensive and, therefore, we approximately invert the two blocks on the main diagonal, and obtain the following block lower triangular preconditioner which we use in our numerical experiments,
\begin{equation} \label{eqn:inexact-block-lower}
	\mathcal{M}_L = 
	\begin{pmatrix}
		Q^{-1}_{\Omega\Omega} &  0  \\
		A_{\Gamma\Omega}  & Q^{-1}_{\Gamma\Gamma}	
	\end{pmatrix}^{-1},
\end{equation}
where $Q_{\Omega\Omega}$ and $Q_{\Gamma\Gamma}$ approximate $\widetilde{S}_{\Omega\Omega}^{-1}$ and $A_{\Gamma\Gamma}^{-1}$, respectively. 
%In practice, we use $\mathcal{M}_L$ as the block preconditioner with proper choices of $Q_{\Omega\Omega}$ and $Q_{\Gamma\Gamma}$.

%\textcolor{red}{TBA: Discussions on other possible choices of Schur complement approximation -- Xiaozhe}

\subsection{Practical Implementation}
In practice, given a right hand side $r = (r_{\Omega}, r_{\Gamma})^T$, an algorithm for the action of the precondtioner $\mathcal{M}_L$ is shown in Algorithm~\ref{alg:action_ML}.
\begin{algorithm}
\caption{Action of the precontioner $\mathcal{M}_L$: $z \gets \mathcal{M}_L r$} \label{alg:action_ML}
\begin{algorithmic}[1]
\State Solve in the domain: $z_{\Omega} \gets Q_{\Omega\Omega} r_{\Omega}$
\State Update the interface residual: $r_{\Gamma} \gets  r_{\Gamma} - A_{\Gamma \Omega} z_{\Omega}$ 
\State Solve on the interface: $z_{\Gamma} \gets Q_{\Gamma \Gamma} r_{\Gamma}$ 
\State Output the update: $z \gets (z_{\Omega}, z_{\Gamma})^T$
\end{algorithmic}
\end{algorithm}
The algorithm requires properly choosen $Q_{\Omega\Omega}$ and $Q_{\Gamma\Gamma}$. 
Since we are solving the scalar elliptic problem, we do this by applying one V-cycle of an algebraic multigrid (AMG) method to $\widetilde{S}_{\Omega\Omega}$ and $A_{\Gamma\Gamma}$, respectively. 
Specifically, our implementation uses smoothed aggregation AMG (SA-AMG) methods to balance the computational complexity and convergence behavior~\cite{vanvek1996algebraic,brezina2001convergence,brezina2012improved,hu2019modifying}. In general, other variants of AMG methods can be applied and similarly the number of V-cycle steps can be modified. 
However, in our numerical experiments, it seems that one V-cycle is sufficient to provide a good approximation and, as a result, leads to an effective preconditioner. 

We comment that $\widetilde{S}_{\Omega\Omega}$ and $A_{\Gamma\Gamma}$ have block structures themselves since we put subdomains $\Omega_i$ and interfaces $\Gamma_i$ of different dimensions together. Therefore, it is possible to design special geometric and algebraic MG methods for approximating $\widetilde{S}_{\Omega\Omega}$ and $A_{\Gamma\Gamma}$. For example, in the setup phase of SA-AMG, one could carefully design coarsening strategies so that aggregations will be constructed within the subdomains and interfaces of the same dimension and then form aggregations that possibly cross different dimensions. In addition, if the subdomains $\Omega_i$ and $\Gamma_i$ are in the one-dimensional space, the cost of directly inverting the corresponding block is relatively negligible. In this work, our choice, SA-AMG implemented in the HAZmath package~\cite{HAZmath}, already provides a good performance, and the specially tailored strategies suggested above do not appear to be necessary. 

\section{Numerical Results}
In this section, we present numerical results to demonstrate the effectiveness of the proposed block preconditioner for solving the linear systems of equations after discretizing the mixed-dimensional scalar elliptic problem~(\ref{eqn:scalar-elliptic-domain}) and~(\ref{eqn:scalar-elliptic-interface}). We use the preconditioner $\mathcal{M}_L$~\eqref{eqn:inexact-block-lower} to accelerate the GMRes method. One V-cycle SA-AMG with one step of Gauss-Seidel iteration as both pre-and post-smoothing steps is used to define ~$Q_{\Omega\Omega}$ and $Q_{\Gamma\Gamma}$, respectively. In all our numerical experiments, we use a zero initial guess, and the GMRes method terminates when the relative residual is smaller than $10^{-6}$. Numerical experiments are conducted on a workstation with an 8-core 3GHz Intel Xeon “Sandy Bridge” CPU and 256 GB of RAM. The software packages used are PorePy~\cite{porepy} (for the discretization of the mixed-dimensional scalar elliptic problems) and HAZmath~\cite{HAZmath} (for the preconditioners and iterative solvers for solving the linear systems of equations). The meshes are generated by Gmsh \cite{geuzaine2009gmsh}, using PorePy's interface to Gmsh to control the mesh size.
The runscripts used to produce the results presented below are available at \cite{runscripts}.

\begin{figure}
    \centering
    \begin{subfigure}[b]{0.32\textwidth}
    \centering
    \includegraphics[scale=0.45]{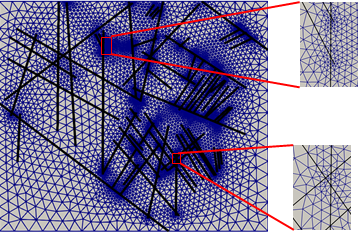}
    \end{subfigure}
    \begin{subfigure}[b]{0.32\textwidth}
    \centering
    \includegraphics[scale=0.3]{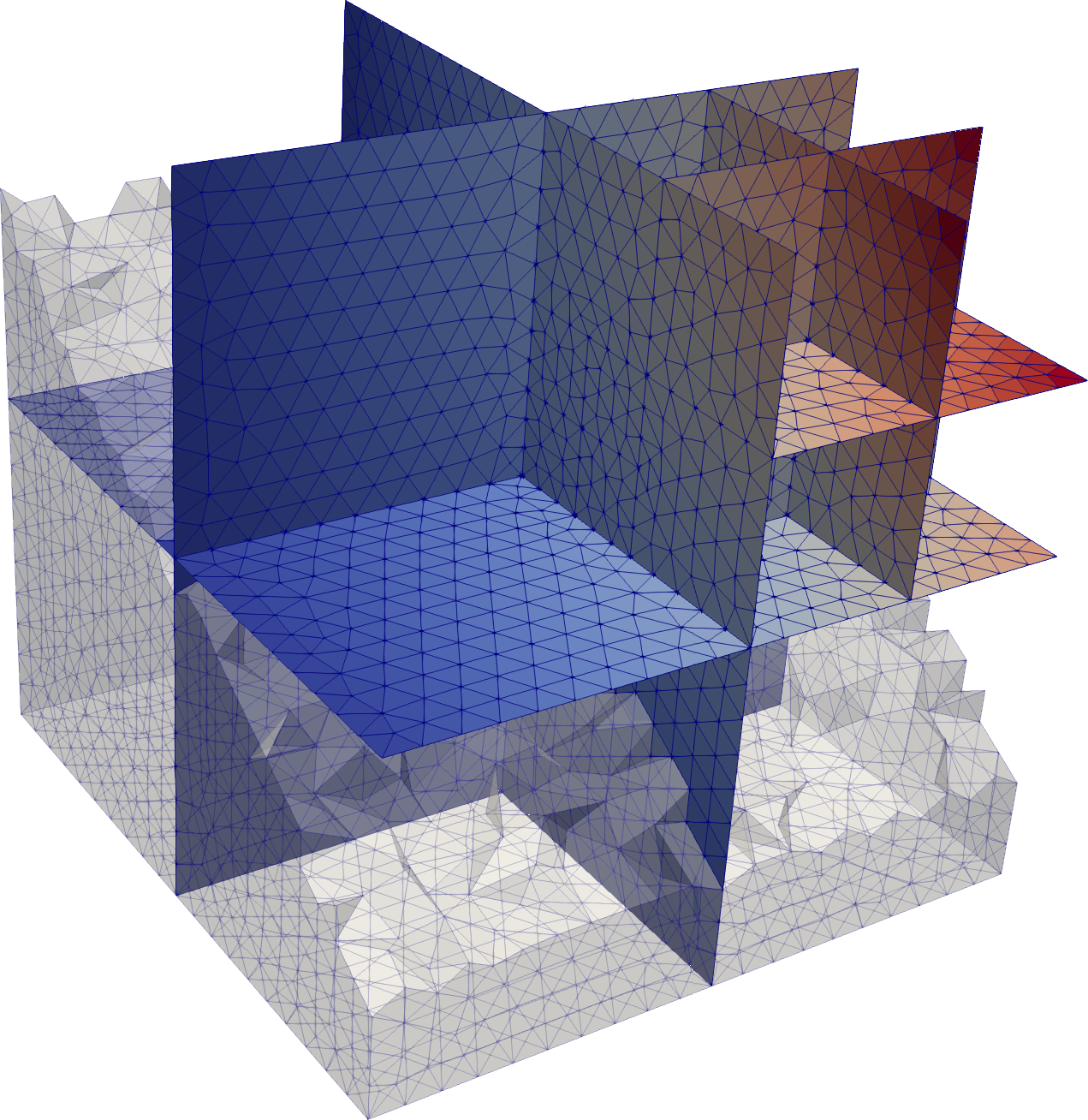}
    \end{subfigure}
    \begin{subfigure}[b]{0.32\textwidth}
    \centering
    \includegraphics[scale=0.25]{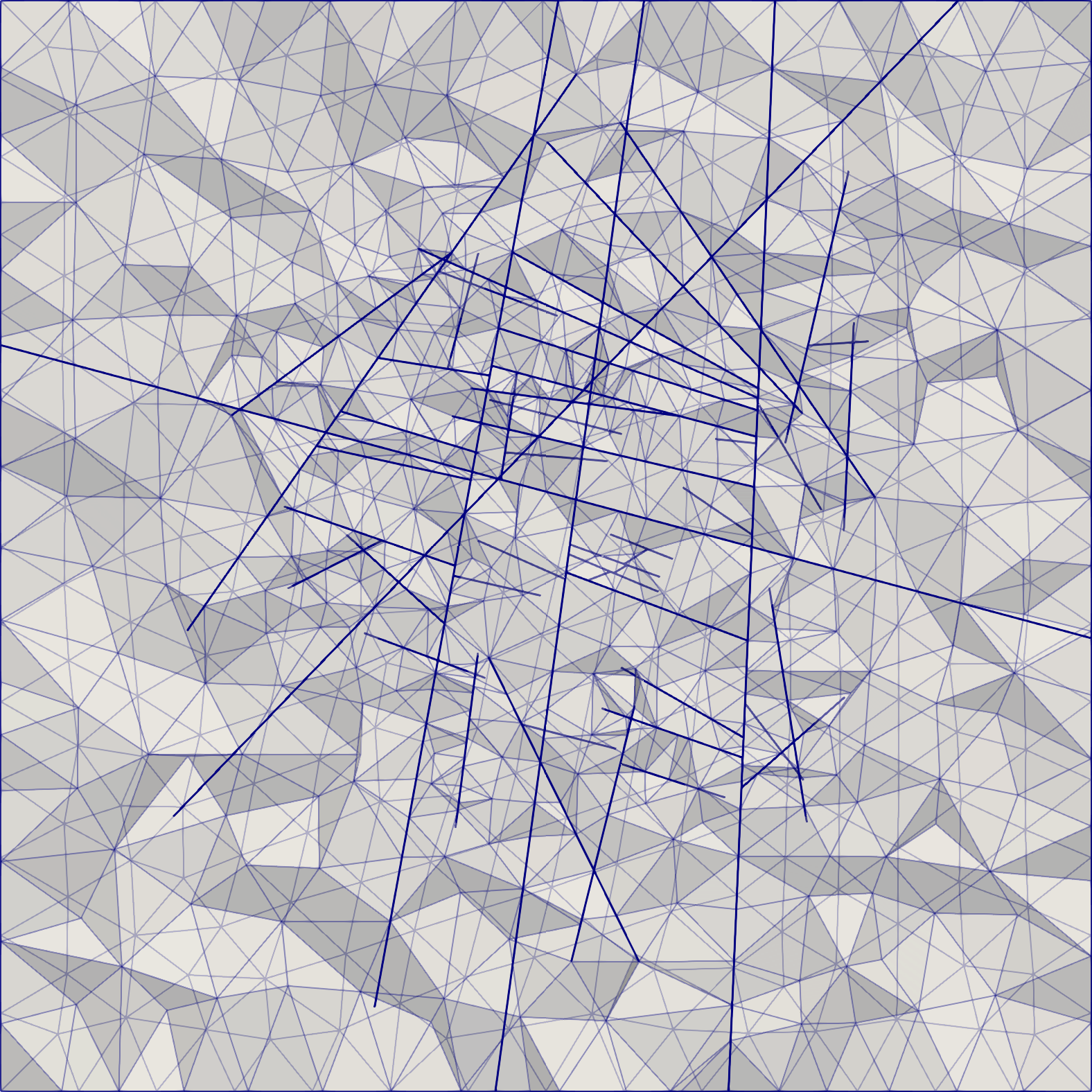}
    \end{subfigure}
    \caption{Left: Fracture geometry (black lines) and grid for the complex 2d case. Middle: Fracture geometry and grid for the regular 3d case; 2d objects are colored by the pressure solution. Right: Fracture network (shown as bold lines) and grid for the 3d field case, top view.}
    \label{fig:2d_geometry}
\end{figure}

\begin{table}[h!]
	\center
	{\small
	\begin{tabular}{|c|c|c|c|c||c|c|c||c|c|c|}
		\hline \hline 
		\multicolumn{2}{|c||}{} &\multicolumn{3}{|c||}{2D Complex} &
		\multicolumn{3}{|c||}{3D Regular} &
		\multicolumn{3}{|c|}{3D Field} \\ \hline 
		$K_{\parallel}$ & $\kappa$  & $h$=20 & $h$=10 & $h$=5 & $h$=0.2 & $h$=0.1 & $h$=0.05 & $h_{f}$=80   & $h_{f}$=40   &  $h_{f}$=20 \\
		\hline \hline
		$10^{-4}$ & $10^{-4}$ & 12  & 12 & 18 & 7  & 7  & 8  & 18  & 16  & 34 \\
		$10^{-4}$ & $1$       & 13  & 12 & 14 & 19 & 20 & 25 & 36  & 51  & 44 \\
		$10^{-4}$ & $10^{4}$  & 16  & 16 & 19 & 14 & 14 & 16 & 68  & 68  & 44 \\
		$1$       & $10^{-4}$ & 10  & 10 & 17 & 8  & 10 & 11 & 13  & 15  & 32 \\
		$1$       & $1$ 	  & 13  & 12 & 13 & 8  & 7  & 8  & 43  & 64  & 50 \\
		$1$       & $10^{4}$  & 16  & 16 & 19 & 13 & 13 & 13 & 62  & 73  & 47 \\
		$10^{4}$  & $10^{-4}$ & 8   & 7  & 10 & 8  & 8  & 8  & 12  & 17  & 28 \\
		$10^{4}$  & $1$       & 28  & 28 & 28 & 8  & 10 & 10 & 32  & 40  & 45  \\
		$10^{4}$  & $10^{4}$  & 23  & 23 & 27 & 10 & 9  & 9  & 181 & 74  & 55\\
		\hline \hline
	\end{tabular}
}
	\caption{Number of iterations of GMRes with block preconditioner $\mathcal{M}_L$ when varying mesh size $h$ (or mesh size of the fracture $h_f$), tangential fracture permeability $K_{\parallel}$, and normal fracture permeability $\kappa$. Left: 2D complex fracture example. Middle: 3D Regular fracture example. Right: 3D field example.
	} \label{tab:all}
\end{table}

\subsection{2D Complex Fracture Example}
For the first example, we choose a 2D example with a complex fracture configuration~\cite{Flemisch2016a} to demonstrate the robustness of the block preconditioner on a realistic fracture network. Such a complex fracture structure often occurs in geological rock simulations, where the geometrical and physical properties of the fracture network can significantly influence the stability of the linear solvers. In particular, as we can see from Figure~\ref{fig:2d_geometry} (Left), the fractured porous medium domain where tips and very acute intersections may decrease the shape regularity of the mesh. Hence, this test case is well suited to show the robustness of the block preconditioner with respect to challenging geometric configurations that are common in realistic fracture geometries.

In the numerical experiments, we set the permeability of the matrix to be $1$. The tangential and normal permeability of the fractures, denoted $K_{\parallel}$ and $\kappa$, respectively, are constants throughout the whole network. To show the robustness of the block preconditioner with respect to the discretization parameter $h$ and physical parameters $K_{\parallel}$ and $\kappa$, we perform a set of tests in which we vary the values of those parameters.

%\begin{table}[h!]
%	\center
%	\begin{tabular}{|c|c|c|c|c|c|}
%		\hline \hline 
%		%\multicolumn{6}{|c|}{2D benchmark complex} \\ \hline 
%		$K_{\parallel}$ & $\kappa$  & $h=40$ & $h=20$ & $h=10$ & $h=5$  \\
%		\hline \hline
%		$10^{-4}$ & $10^{-4}$ & 12  & 12  & 12 & 18   \\
%		$10^{-4}$ & $1$       & 14  & 13  & 12 & 14   \\
%		$10^{-4}$ & $10^{4}$  & 15  & 16  & 16 & 19   \\
%		$1$       & $10^{-4}$ & 10  & 10  & 10 & 17   \\
%		$1$       & $1$ 	  & 13  & 13  & 12 & 13   \\
%		$1$       & $10^{4}$  & 15  & 16  & 16 & 19   \\
%		$10^{4}$  & $10^{-4}$ & 8   & 8   & 7  & 10   \\
%		$10^{4}$  & $1$       & 27  & 28  & 28 & 28   \\
%		$10^{4}$  & $10^{4}$  & 21  & 23  & 23 & 27   \\
%		\hline \hline
%	\end{tabular}
%\caption{2D Complex fracture example: number of iterations of GMRes with block preconditioner $\mathcal{M}_L$ when varying mesh size $h$, tangential fracture permeability $K_{\parallel}$, and normal fracture permeability $\kappa$.} \label{tab:2D-complex-fracture}
%\end{table}

Table~\ref{tab:all}(Left) summarizes the numerical results. Each row in the table represents a set of tests where the tangential and normal fracture permeability values are fixed and the mesh size $h$ is changing. 
It should be mentioned that these are mesh sizes prescribed to PorePy, but in practice, local geometric details combined with Gmsh's internal functionality to generate high-quality meshes may lead to further local grid refinement.
As we can see, the numbers of GMRes iterations grow slightly but remain under $ 30$ iterations, which demonstrates the robustness of the preconditioner with respect to the mesh size $h$. On the other hand, each column of the table presents a set of tests where the mesh size is fixed, and the tangential and normal fracture permeability values vary. Again, the numbers of GMRes iteration remain stable as expected, which demonstrate the robustness of the block preconditioner with respect to the physical parameters.   Overall, Table~\ref{tab:all}(Left) illustrates the robustness of the block preconditioner for this 2D complex fracture benchmark.

\subsection{3D Regular Fracture Example}
Now we consider a 3D problem from a benchmark study~\cite{berre2021verification}. The 3D geometry is a unit cube as shown in Figure~\ref{fig:2d_geometry} (Middle). The fracture network consists of $9$ fracture planes, $69$ intersection lines, and $27$ intersection points. Again, to show the robustness of the block preconditioner with respect to the discretization and physical parameters, we vary the values of the tangential and normal fracture permeability and the mesh size, keeping the permeability of the porous medium unitary.

%\begin{table}[h!]
%	\center
%	\begin{tabular}{|c|c|c|c|c|}
%		\hline \hline 
%		%\multicolumn{5}{|c|}{3d regular} \\ \hline 
%		$K_{\parallel}$ & $\kappa$  & $h=0.2$ & $h=0.1$ & $h=0.05$   \\
%		\hline \hline
%		$10^{-4}$ & $10^{-4}$ & 7  & 7  & 8    \\
%		$10^{-4}$ & $1$       & 19 & 20 & 25   \\
%		$10^{-4}$ & $10^{4}$  & 14 & 14 & 16   \\
%		$1$       & $10^{-4}$ & 8  & 10 & 11   \\
%		$1$       & $1$ 	  & 8  & 7  & 8    \\
%		$1$       & $10^{4}$  & 13 & 13 & 13   \\
%		$10^{4}$  & $10^{-4}$ & 8  & 8  & 8    \\
%		$10^{4}$  & $1$       & 8  & 10 & 10   \\
%		$10^{4}$  & $10^{4}$  & 10 & 9  & 9    \\
%		\hline \hline
%	\end{tabular}
%\caption{3D regular fracture example: number of iterations of FGMRes with block preconditioner $\mathcal{M}_L$ when varying mesh size $h$, tangential fracture permeability $K_{\parallel}$, and normal fracture permeability $\kappa$.} \label{tab:3D-regular-fracture}
%\end{table}

Table~\ref{tab:all}(Middle) shows the number of GMRes iterations for this 3D regular fracture example.  We can see that the numerical results are consistent with our 2D example, i.e., the block preconditioners show robustness with respect to the mesh size $h$, the tangential fracture permeability $K_{\parallel}$, and normal fracture permeability $\kappa$.   

\subsection{3D Field Example}
Our last example is a simulation of a 3D field benchmark with a realistic fracture network~\cite{berre2021verification}. The domain is $\Omega = (-500, 350) \times (100,1500) \times (-100, 500)$, and a cross-section of the domain is shown in Figure~\ref{fig:2d_geometry} (Right) where we can see the complex fracture network and computational grid. In this example, we again set the permeability in the matrix to be $1$ and vary the values of the tangential and normal fracture permeability and the mesh size in our numerical experiments, which are the same as in previous examples. 
However, due to the complexity of this example and computational cost consideration, we fix the far-field mesh size of the matrix to be $100$ and vary the degree of mesh refinement towards the fractures, mesh size in the fractures (denoted by $h_{f}$) only. 
Larger values of $h_{f}$ give more irregular the elements, especially near the intersection of the fractures and many tightly packed fractures, thus, the performance of the linear solver is expected to improve as $h_f$ is decreased.

%\begin{table}[h!]
%	\center
%	\begin{tabular}{|c|c|c|c|c|c|}
%		\hline \hline 
%		%\multicolumn{6}{|c|}{3d field ($h=100$)} \\ \hline 
%		$K_{\parallel}$ & $\kappa$  & $h_{frac}=100$  &  $h_{frac}=80$   & $h_{frac}=40$   &  $h_{frac}=20$  \\   
%		%\hline 
%		%&  & \#DoF = $51,558$  & \#DoF = $61,330$  &  \#DoF = $149,569$ &  \#DoF = $563,435$  \\
%		\hline \hline        
%		$10^{-4}$ & $10^{-4}$ & 15  & 18  & 16  & 34 \\ 
%		$10^{-4}$ & $1$       & 47  & 36  & 51  & 44  \\
%		$10^{-4}$ & $10^{4}$  & 73  & 68  & 68  & 44 \\
%		$1$  	  & $10^{-4}$ & 13  & 13  & 15  & 32 \\
%		$1$       & $1$       & 55  & 43  & 64  & 50 \\
%		$1$       & $10^{4}$  & 85  & 62  & 73  & 47 \\ 
%		$10^{4}$  & $10^{-4}$ & 12  & 12  & 17  & 28 \\
%		$10^{4}$  & $1$       & 34  & 32  & 40  & 45 \\
%		$10^{4}$  & $10^{4}$ & 193 & 181  & 74  & 55  \\
%		\hline \hline
%	\end{tabular}
%\caption{3D field example: number of iterations of FGMRes with block preconditioner $\mathcal{M}_L$ when varying mesh size $h$, tangential fracture permeability $K_{\parallel}$, and normal fracture permeability $\kappa$.} \label{tab:3D-field-fracture}
%\end{table} 

The results for the 3D field example are reported in Table~\ref{tab:all}(Right).  As expected, when $h_{f}$ is large, the mesh quality influences the overall performance, and the number of GMRes iterations varies quite a bit, see the  columns of $h_f=80$ in Table~\ref{tab:all}(Right).  When the fracture mesh size $h_{f}$ gets smaller, the number of GMRes iterations stabilizes, and the performance is robust with respect to the discretization and physical parameters as we expected.  Due to the complex geometry, the overall number of iterations is higher than the 3D regular fracture example.  Therefore, for complex fracture networks, we suggest investing in constructing a more regular mesh of the fractured porous medium and then applying the proposed block preconditioners in the iterative solvers.

\section{Conclusions}
Based on the block structure of the linear systems arising from discretizing the mixed-dimensional scalar elliptic problems, we are able to develop block preconditioners based on approximate factorization. We first properly approximate the Schur complement to obtain a block preconditioner and then apply the AMG methods to invert the diagonal blocks in our practical implementation. Several benchmarks in 2D and 3D are considered. From the numerical results, the GMRes methods accelerated by our proposed preconditioner are robust with respect to the physical and discretization parameters and complex fracture structures, making it attractive for real-world applications. 

\section*{Acknowledgments}
This work was financed in part by Norwegian Research Council grant number 308733.

\printbibliography

\end{document}